\documentclass[10pt,twoside]{article}
\usepackage{amssymb}
\usepackage{amsmath} 
\usepackage{amsfonts} 
\usepackage{amssymb} 
\usepackage{latexsym}
\newtheorem{theorem}{Theorem}[section]

\newtheorem{lemma}{Lemma}[section]
\newtheorem{proposition}{Proposition}[section]

\newtheorem{corollary}{Corollary}[section]
\def\proof{\mbox {\it Proof.~}}

\renewcommand{\a }{\alpha } 
\renewcommand{\b }{\beta } 
\renewcommand{\d}{\delta }

\newcommand{\D }{\Delta }
 
\newcommand{\e }{\varepsilon }

\renewcommand{\l }{\lambda }

\newcommand{\n }{\nabla }

\newcommand{\be}{\begin{equation}} 
\newcommand{\ee}{\end{equation}} 
 
\newenvironment{pfn}[1]{\noindent{\bf Proof of {#1}\enspace}}{
\hfill$\Box$\medskip} 
\newcommand{\R}{\mathbb{R}} 
 
\newcommand{\N}{\mathbb{N}}


\begin{document}
\title{
\vspace{0.5in}
{\bf\Large Prescribing the Scalar Curvature under Minimal\\
 Boundary Conditions on the Half Sphere\thanks{ The main part of this
 work  was done while the authors enjoyed the hospitality of Mathematisches Forschungsinstitut at Oberwolfach. Their stay was supported by Volkswagen-Stiftung (RIP-Program at Oberwolfach)} }}
\author{{\bf\large Mohamed Ben Ayed }\hspace{2mm}\\
{\it\small Department de Mathematiques},\\ {\it\small Faculte des Sciences de Sfax}, {\it\small Route Soukra, Sfax, Tunisia}\\
{\it\small e-mail: Mohamed.Benayed@fst.rnu.tn}\vspace{1mm}\\
{\bf\large Khalil El Mehdi }\hspace{2mm}\\
{\it\small Centre de Math{\'e}matique },\\ {\it\small Ecole Polytechnique}, {\it\small 91128 Palaiseau, France}\\
{\it\small e-mail: elmehdi@math.polytechnique.fr}\vspace{1mm}\\
{\it\small and }\\
{\it\small Faculte des Sciences et techniques  },\\ {\it\small Universite de Nouakchott, Mauritanie} \\
{\bf\large Mohameden Ould Ahmedou}\vspace{1mm}\\
{\it\small Mathematisches Institut},\\  {\it\small Beringstrasse 4},\\
{\it\small D-53115 Bonn, Germany }\\
{\it\small e-mail: ahmedou@math.uni-bonn.de}\vspace{1mm}}

\date{ }
\maketitle
\begin{center}
{\bf\small Abstract }

\vspace{3mm}
\hspace{.05in}\parbox{4.5in}
{{\small   This paper is devoted to the problem of prescribing the scalar curvature under zero boundary conditions. Using dynamical and topological methods involving the study of critical points at infinity of the associated variational problem, we prove some existence results on the standard half sphere                               }}
\end{center}


\section{Introduction and the Main Results }
In this paper we study some nonlinear problem arising from conformal geometry.
Precisely, consider a Riemannian manifold with boundary $(M^n,g)$ of dimension 
$n\geq 3$ and take $\tilde{g}=u^{4/(n-2)}g$, be a conformal metric to $g$, 
where $u$ is a smooth positive function, then the following equations relate 
the scalar curvatures $R_g$, $R_{\tilde{g}}$ and the mean curvatures of the 
boundary $h_g$, $h_{\tilde{g}}$, with respect to $g$ and $\tilde{g}$ 
respectively. 
\begin{eqnarray*}
(P_1)\quad \left\{
\begin{array}{ccccc} 
-c_n\D _gu+R_gu&=&R_{\tilde{g}}u^{\frac{n+2}{n-2}}&\mbox{ in }& M\\
\frac{2}{n-2}\frac{\partial u}{\partial \nu}+h_g u &=&
h_{\tilde{g}}u^\frac{n}{n-2}&
\mbox{ on }& \partial M
\end{array}
\right.
\end{eqnarray*}
where $c_n=4(n-1)/(n-2)$ and $\nu$ denotes the outward normal vector with 
respect to the metric $g$.\\
In view of the above equations, a natural question is whether it is possible to prescribe both the scalar 
curvature and the boundary mean curvature, that is : given two functions 
$K: M\to \R$ and $H: \partial M\to \R$, does exists a metric $\tilde{g}$
conformally equivalent to $g$ such that $R_{\tilde{g}}=K$ and 
$h_{\tilde{g}}=H$? \\
According to equations $ (P_1)$, the problem is equivalent to finding a smooth 
positive solution $u$ of the following equation
\begin{eqnarray*}
(P_2)\quad \left\{
\begin{array}{ccccc} 
-c_n\D _gu+R_gu&=& K u^{\frac{n+2}{n-2}}&\mbox{ in }& M \\
\frac{2}{n-2}\frac{\partial u}{\partial \nu}+h_g u &= & 
 H u^\frac{n}{n-2}&
\mbox{ on }& \partial M.
\end{array}
\right.
\end{eqnarray*}
Such a problem was studied in \cite{ALM} \cite{C},\cite{DMA}  \cite{E1}, \cite{E2},  \cite{HL1}, \cite{HL2} \cite{L1} . Yanyan Li \cite{L1}, and Djadli-Malchiodi-Ould Ahmedou  \cite{DMA} 
studied this  problem  when the manifold is the three dimensional standard half
sphere. Their approach involves a fine blow up analysis of some subcritical 
approximations and the use of the topological degree tools.

Regarding the above  problem it is well known that the most interesting case is the so called positive one, that is when the quadratic part of the associated Euler functional  is positive definite. Another interesting case is when a noncompact group of conformal transformations acts on the equation leading to topological obstructions. The half sphere represents the simplest case where such a noncompactness occurs, and in this paper we consider the case of the standard half sphere under minimal boundary conditions:

More precisely, let 
$$
S_+^n=\{x\in \R^{n+1}\, /\, |x|=1,\, x_{n+1}>0\}
$$
$n\geq 3$ . 
Given a $C^2$ function $K$ on 
$\overline{S_+^n}$, we look for conditions on $K$ to ensure the existence of a  positive solution of the problem
\begin{eqnarray*}
(1)\quad \left\{
\begin{array}{ccccc} 
-\D _gu+ \frac{n(n-2)}{4}u&=& K u^{\frac{n+2}{n-2}}&\mbox{ in }& S_+^n\\
\frac{\partial u}{\partial \nu}&=&0& 
\mbox{ on }& \partial S_+^n
\end{array}
\right.
\end{eqnarray*}
where $g$ is the standard metric of $S_+^n$.\\
Problem (1) is in some sense related to the well known scalar curvature 
problem on $S^n$
\begin{eqnarray*}
(2)\qquad-\D _gu+ \frac{n(n-2)}{4}u= K u^{\frac{n+2}{n-2}}\quad \mbox{ in } S^n
\end{eqnarray*}
to which much work has been devoted. For details please see 
  \cite{B1},\cite{BC1} , \cite{BCCH},\cite{BCH},   \cite{CY},    \cite{CGY},  \cite{H}, \cite{L2}, \cite{L3},\cite{SZ} 
  and the references therein.\\
As for (2), there are topological obstructions of Kazdan-Warner type to solve 
(1) (see  \cite{BP}) and so a naturel question arises: under which 
conditions on $K$, (1) has a positive solution. We propose to handle such a
 question, using some topological and dynamical tools of the theory of 
critical points at infinity, see Bahri \cite{B1} , \cite{B2}.\\
Our approch follows closely the ideas developped in Aubin-Bahri \cite{AB}, 
Bahri \cite{B1} and Ben Ayed-Chtioui-Hammami \cite{BCH} where the problem of 
prescribing the scalar curvature on closed manifolds was studied using 
some algebraic topological tools. The main idea is to use the difference of 
topology between the level sets of the function $K$ to produce a critical 
point of the Euler functional $J$ associated to (1) and the main issue is 
under which conditions on $K$, a topological accident between the level sets 
of $K$ induces a topological accident between the level sets of $J$. Such an accident is sufficient to prove the existence of a critical point
when some compactness conditions are satisfied. However our problem presents a 
lack of compactness due to the presence of critical points at infinity, that 
is noncompact orbits for the gradient of $J$ along which $J$ is bounded and 
its gradient goes to zero. Therefore a careful study of such noncompact orbits 
is necessary, in order to take into account their contribution to the 
difference of topology between the level sets of $J$.\\
In order to state our main results, we need to introduce the assumptions 
that we are using in our results.\\
${\bf (A_1)}$\hskip 0.3cm We assume  that 
$K_1=K_{|\partial S_+^n}$ has only nondegenerate critical points $y_0,...,y_s$, where 
$y_0$ is the absolute maximum, such that 
$$
K(y_0)\geq K(y_1)\geq ...\geq K(y_l)> K(y_{l+1}\geq ...\geq K(y_s)
$$
with
$$
\frac{\partial K}{\partial \nu}(y_i)>0,\, \mbox{ for } 0\leq i\leq l, \qquad 
\frac{\partial K}{\partial \nu}(y_i)\leq 0, \,  \mbox{ for } l+1\leq i\leq s.
$$ 
${\bf (A_2)}$\hskip 0.3cm Assume that there exists 
$c$  a postive constant such that $c < K(y_l)$, and  every $y$  critical point of $K$, $K(y) < c $.\\
${\bf (A_3)}$ \hskip 0.3cm Let $Z$ be a pseudogradient of $K_1$, of Morse-Smale 
type (that is the intersections of the stable and the unstable manifolds of the critical points of $K_1$ are transverse.)\\
Set 
$$
X=\overline{\cup_{0\leq i\leq l}W_s(y_i)}
$$ 
where $W_s(y_i)$ is the stable manifold of $y_i$ for $Z$.\\
We assume that $X$ is not contractible and denote by $m$ the dimension of the 
first nontrivial reduced homological group.\\
${\bf (A_4)}$ \hskip 0.3cm Assume that $X$ is contractible in 
$K^c=\{x\in S_+^n\, /\, K(x)\geq c\}$, where $c$ is defined in the 
assumption $(A_2)$.\\
Now, we are able to state our first main results.
\begin{theorem}\label{t:11}
Assume that $n\geq 4$. Then, under the assumptions $(A_1)$, $(A_2)$, $(A_3)$, $(A_4)$, there exists a 
constant $c_0$ independent of $K$ such that if $K(y_0)/c\leq 1+c_0$, 
then (1) has a solution.
\end{theorem}
\begin{corollary}\label{c:12}
The solution obtained in Theorem \ref{t:11} has an augmented Morse \\
index $\geq m$.
\end{corollary}

Next, we state another kind of existence results for problem (1) based on a `` topological invariant'' for some Yamabe type problems  introduced by Bahri see \cite{B1}. \\
To state these results, we need to introduce the assumptions that we will be
using  and some notations.\\
${\bf (H_1)}$\hskip 0.3cm We assume that $K_1$ has only nondegenerate
critical points and we assume that there exists $y_0 \in \partial
S^n_+$ such that $y_0$ is the absolute maximum of $K_1$ and $(\partial K/\partial\nu )(y_0) > 0$.\\
${\bf (H_2)}$\hskip 0.3cm $W_s (y_i)\cap W_u(y_j) = \emptyset $ for
any $i$ such that $(\partial K/\partial\nu )(y_i) > 0$ and for any $j$
such that $(\partial K/\partial\nu )(y_j) < 0$.\\
For $ k \in \{1, 2, ..., n-1 \} $, we define $X$ as 
$$
X=\overline{W_s(y_{i_0})}
$$
where $y_{i_0}$ satisfies 
$$
K_1(y_{i_0})= \max \{K_1(y_i), \, /\, \mbox{ind}(y_i)=n-1-k, \, 
({\partial K}/{\partial \nu})(y_i)> 0\}
$$ 
(Here ind$(y_i)$ denotes the Morse index of $y_i$ for the function $K_1$).\\
${\bf (H_3)}$\hskip 0.3cm We assume that $X$ is without boundary. \\
We denote by $C_{y_0}(X)$ the following set 
$$
C_{y_0}(X)=\{\a \d _{y_0}+(1-\a )\d _x \, / \, \a \in [0,1],\, x\in X \}.
$$
For $\l $ large enough, we introduce a map $f_\l : C_{y_0}(X)\to \Sigma ^+$,
defined by
$$
C_{y_0}(X) \ni (\a \d _{y_0}+(1-\a )\d _x)  \longrightarrow 
\frac{\a \d _{(y_0,\l )}
+(1-\a )\d _{(x,\l )}}{|\a \d _{(y_0,\l )}+(1-\a )\d _{(x,\l )}|}\in
\Sigma ^+.
$$
Then $ C_{y_0}(X)$ and $f_\l ( C_{y_0}(X))$ are manifolds in dimension $k+1$, 
that is, their singularities arise in dimension $k-1$ and lower, see \cite{B1}. Observe that 
 $ C_{y_0}(X)$ and $f_\l ( C_{y_0}(X))$ are contractible while $X$ is 
not contractible.\\
For $\l $ large enough, we also define the intersection number(modulo 2) of 
$f_\l (C_{y_0}(X))$ with $ W_s(y_0,y_{i_0})_\infty$
$$
\mu (y_0)=f_\l (C_{y_0}(X)). W_s(y_0,y_{i_0})_\infty
$$ 
where $W_s(y_0,y_{i_0})_\infty$ is the stable manifold of 
$(y_0,y_{i_0})_\infty$ for a 
decreasing pseudogradient $V$ for $J$ which is transverse to 
$f_\l (C_{y_0}(X))$. Thus this number is well defined \cite{M}.\\
We then have the following result:
\begin{theorem}\label{t:13}
Assume that $n\geq 4$. Under assumptions $(H_1)$, $(H_2)$ and $(H_3)$,
if $\mu (y_0)=0$ then (1) has a solution of index $k$ or $k+1$.
\end{theorem}
Now, we state a statement  more general than Theorem \ref{t:13}.  For this we
define $X$ to be 
$$
X=\overline{\cup_{y_i\in B_k} W_s(y_i)}, \quad \mbox{ with } B_k=\{y_i \, /
\, \mbox{ind}(y_i)=n-1-k \mbox{ and } ({\partial K}/{\partial \nu})(y_i)> 0\}.
$$
For $y_i\in B_k$, we define, for $\l $ large enough, the intersection
number(modulo 2) 
$$
\mu_i (y_0)=f_\l (C_{y_0}(X)). W_s(y_0,y_i)_\infty.
$$
By the above arguments, this number is well defined \cite{M}. \\
Then we have the following theorem
\begin{theorem}\label{t:14}
Assume that $n\geq 4$. Under assumptions $(H_1)$, $(H_2)$ and $(H_3)$,
if $\mu_i=0$ for each $y_i\in B_k$, then (1) has a solution of index $k$ 
or $k+1$.
\end{theorem}
The remainder of the present paper is organized as follows. In section 2, we set up 
the variational structure and recall some preliminaries. In section 3, we 
perform an expansion of the Euler functional associated to
 (1) and its gradient near the potential  critical 
points at infinity, then we prove a Morse Lemma at infinity in section 4. 
In section 5, we provide the proof of Theorem \ref{t:11} and Corollary \ref{c:12}, while section 6 is devoted to the proof of Theorems 
\ref{t:13} and \ref{t:14}.

\section { Variational Structure and Preliminaries }

In this section we recall the functional setting and the variational problem 
and its main features. Problem (1) has a variational structure. The 
functional is 
$$
J(u)= \frac {\int_{S_+^n}|\n u|^2+\frac {n(n-2)}{4}\int_{S_+^n}u^2}
{\left(\int_{S_+^n}Ku^{\frac{2n}{n-2}}\right)^{\frac {n-2}{n}}}
$$
defined on $H^1(\overline{S_+^n},\R)$ equiped with the norm 
$$
||u||^2=\int_{S_+^n}|\n u|^2+ \frac{n(n-2)}{4} \int_{S_+^n}u^2.
$$
We denote by $\Sigma$ the unit sphere of $H^1(\overline{S_+^n},\R)$ and we set
$\Sigma ^+=\{u \in \Sigma \, / \, u\geq 0\}$.\\
The Palais-Smale condition fails to be satisfied for $J$ on $\Sigma ^+$. Its failure has been 
studied by various authors (see Brezis-Coron \cite{BrC}, Lions \cite{L}, Struwe \cite{S}).\\
In order to characterize the sequences failing the Palais-Smale condition, we 
need to introduce some notations.\\
For $a\in \overline{S_+^n}$ and $\l >0$, let 
$$
\tilde\d _{a,\l }(x)= c_0\frac {\l ^{\frac{n-2}{2}}}{(\l ^2+1+(\l ^2-1)
\cos d(a,x))^{\frac{n-2}{2}}}
$$
where $d$ is the geodesic distance on $(S_+^n,g)$ and $c_0$ is chosen so that 
$$
-\D \tilde\d _{a,\l } + \frac{n(n-2)}{4}\tilde\d _{a,\l }=
\tilde\d _{a,\l }^{\frac{n+2}{n-2}}, \quad \mbox{ in } S_+^n.
$$
For $\e >0$ and $p\in \N^*$, let us define  
\begin{align*}
V(p,\e )=& \{u\in \Sigma ^+/\exists \, a_1,...,a_p \in \overline{S_+^n}, \exists \, 
\l _1,...,\l _p >0, \exists \, \a _1,...,\a _p>0 \\
 & \mbox{ s.t. } ||u-\sum_{i=1}^p\a _i\tilde\d _i||<\e \mbox{ and } |\frac 
{\a _i ^{\frac{4}{n-2}}K(a_i)}{\a _j ^{\frac{4}{n-2}}K(a_j)}-1|<\e ,\\ 
 & \l _i>\e ^{-1}, \e _{ij}<\e \mbox{ and }\l _id_i<\e \mbox{ or }
\l _id_i >\e ^{-1}\}
\end{align*}
where $\tilde\d _i=\tilde\d _{a_i,\l _i}$, $d_i=d(a_i,\partial S_+^n)$ and 
$\e _{ij}=(\l _i/\l _j +\l _j/\l _i + \l _i\l _jd^2(a_i,a_j))^{\frac{2-n}{2}}$.
The failure of Palais-Smale condition can be described as follows:
\begin{proposition}\label{p:21} (see \cite{BC2}, \cite{L} and \cite{S})
Assume that $J$ has no critical point in $\Sigma ^+$ and let $(u_k)\in 
\Sigma ^+$ be a sequence such that $J(u_k)$ is bounded and $\n J(u_k)\to 0$.
Then, there exist an integer $p\in \N^*$, a sequence $\e _k>0$ ($\e _k\to 0$)
and an extracted subsequence of $u_k$, again denoted $(u_k)$, such that 
$u_k\in V(p,\e _k )$.
\end{proposition}
Now, we consider the following subset of $V(p,\e )$
$$
V_b(p,\e )=\{u\in V(p,\e )\, /\, \l _id_i<\e \}.
$$
The following lemma defines a parametrization of the set $V_b(p,\e )$.
\begin{lemma}\label{l:22} (see \cite{B2}, \cite{BC2}, \cite{R})
There is $\e _0>0$ such that if $\e <\e _0$ and $u\in V_b(p,\e )$, then the 
problem 
$$
\min\{||u-\sum_{i=1}^p\a _i\tilde\d _i||,\, \a _i>0,\, \l _i>0,\, a_i\in 
\partial S_+^n\}
$$
has a unique solution (up to permutation). In particular, we can write 
$u\in V_b(p,\e )$ as follows 
$$
u=\sum_{i=1}^p\bar{\a }_i\tilde\d _{\bar{a}_i,\bar{\l }_i}+v,
$$
where $(\bar{\a }_1,...,\bar{\a }_p,\bar{a}_1,...,\bar{a}_p,\bar{\l }_1,...,
\bar{\l }_p)$ is the solution of the minimization problem and where 
$v\in H^1(S_+^n)$ such that for each $i=1,...,p$ 
$$
(v,\tilde\d _i)=(v,\partial \tilde\d _i/\partial \l _i)=(v,
\partial \tilde\d _i/\partial a _i)=0.
$$
Here $(.,.)$ denoted the scalar inner defined on $H^1(\overline{S_+^n})$ by 
$$
(u,v)=\int_{S_+^n}\n u \n v + \frac{n(n-2)}{4}  \int_{S_+^n}uv.
$$
\end{lemma}
Before ending this section, we mention that it will be convenient to
perform some stereographic projection in order to reduce our problem
to $\R^n_+$. Let $D^{1,2}(\R^n_+)$ denote the completion of $C^\infty
_c (\bar{\R^n_+})$ with respect to Dirichlet norm. The stereographic
projection $\pi _a$ throught a point $a \in \partial S^n_+$ induces an
isometry $i : H^1(S^n_+) \to D^{1,2}(\R^n_+) $ according to the
following formula
$$
(iv)(x)= \left(\frac{2}{1+|x|^2}\right)^{(n-2)/2}v(\pi _a^{-1}(x)),
\qquad v\in H^1(S^n_+), \, x\in \R^n.
$$
In particular, one can check that the following holds true, for every
$v\in H^1(S^n_+)$
$$
\int _{S^n_+}(|\n v|^2 + \frac{n(n-2)}{4}v^2) = \int _{\R^n_+}|\n
(iv)|^2 \qquad \mbox{and } \int _{S^n_+}|v|^{\frac{2n}{n-2}}= \int
_{\R^n_+}|iv|^{\frac{2n}{n-2}}.
$$
In the sequel, we will identify the function $K$ and its composition
with the stereographic projection $\pi _a$. We will also identify a
point $b$ of $S^n_+$ and its image by $\pi _a$. These
facts will be assumed as understood in the sequel.

\section{ Expansion of $J$ and its gradient at infinity }

This section is devoted to an useful expansion of $J$ and its gradient  near a potential boundary
critical point at infinity consisting of two masses.   
\begin{proposition}\label{p:31}
For $\e >0$ small enough and 
$u=\sum_{i=1}^2\a _i\tilde\d _{a_i,\l _i}+v\in V_b(2,\e )$, 
we have the following expansion
\begin{align*}
J(u)=& \frac{(\a _1^2+\a _2^2)(S_n/2)^{(2/n)}}{(\a _1^{\frac{2n}{n-2}}K(a_1)+
\a _2^{\frac{2n}{n-2}}K(a_2))^{\frac{n-2}{n}}}\left[1+\frac{4(n-2)}{n\b}c_1
\sum_{i=1}^2\frac{\a _i ^{\frac{2n}{n-2}}}{\l _i}\frac{\partial K}{\partial 
\nu}(a_i)\right.\\
 & -2c_2\a _1\a _2\e _{12}\left(-\frac{1}{\gamma}+\frac{1}{\b }(\sum_{i=1}^2
\a _i ^{\frac{4}{n-2}}K(a_i))\right)+f(v)+Q(v,v)\\
 & \left.+O\left(\e _{12}^{\frac{n}{n-2}}log(\e _{12}^{-1})+\sum(\frac{1}
{\l _i ^2}+\frac{\e _{12}}{\l _i}(log(\e _{12}^{-1}))^{\frac{n-2}{n}})+
||v||^{\inf(3,\frac{2n}{n-2})}\right)\right]
\end{align*} 
where 
\begin{align*}
Q(v,v) = & \frac{1}{\gamma}||v||^2-\frac{n+2}{n-2}\frac{1}{\b} \int_{S_+^n}K(\a _1
\tilde\d _1+\a _2 \tilde\d _2)^{\frac{4}{n-2}}v^2,  \\ 
f(v) = & -\frac{2}{\b }\int_{S_+^n}K(\a _1\tilde\d _1+\a _2
\tilde\d _2)^{\frac{n+2}{n-2}}v,\qquad S_n = c_0^{\frac{2n}{n-2}}
\int_{\R^n}\frac{dx}{(1+|x|^2)^n},\\ 
\b  = &\frac{ S_n}{2}(\a _1^{\frac{2n}{n-2}}K(a_1)+\a _2^{\frac{2n}{n-2}}
K(a_2)), \qquad \gamma = \frac{S_n}{2}(\a _1^2+\a _2^2),\\
 c_1 = & c_0^{\frac{2n}{n-2}}\int_{\R_+^n}\frac{x_ndx}{(1+|x|^2)^n},\qquad  
\qquad \qquad c_2=c_0^{\frac{2n}{n-2}}\int_{\R^n}\frac{dx}{(1+
|x|^2)^{\frac{n+2}{2}}.}
\end{align*}
\end{proposition} 
\proof We need to estimate 
$$
 N(u)= ||u||^2 \mbox{ and } D^{\frac{n}{n-2}}=\int_{S_+^n} K(x) 
u^{\frac{2n}{n-2}}. 
$$
In order to simplify the notations, in the remainder, we write
$\tilde\d _i$ instead of $\tilde\d _{a_i, \l _i}$.\\
We  now have 
$$
N(u)= \a _1^2||\tilde\d _1||^2+\a _2^2||\tilde\d _2||^2+||v||^2+2\a _1
\a _2(\int_{S_+^n}\n \tilde\d _1\n \tilde\d _2 +\frac{n(n-2)}{4}
\int_{S_+^n}\tilde\d _1\tilde\d _2)
$$
Observe that 
\begin{eqnarray}\label{e:o1}
||\tilde\d ||^2=\int_{\R_+^n}|\n \d |^2 = \frac{S_n}{2}
\end{eqnarray}
and 
$$
\int_{S_+^n}\n \tilde\d _1\n \tilde\d _2 +\frac{n(n-2)}{4}
\int_{S_+^n}\tilde\d _1\tilde\d _2= 
\int_{\R_+^n}\n \d _1\n \d _2 =\int_{\R_+^n}\d _1^{\frac{n+2}{n-2}} \d _2
$$
where $\d _i$ denotes $  \d _{a_i,\l _i }$ and, for $a \in \R^n$ and $\l >
0$,   $ \d _{a,\l }$ denotes the family of solutions of Yamabe problem on 
$ \R^n$ defined by 
$$
 \d _{a,\l }(x)=c_0\frac{\l ^{\frac{n-2}{2}}}{(1+\l ^2|x-a|^2)^{\frac{n-2}{2}}}.
$$
A computation similar to the one performed in \cite{B2}) shows that 
\begin{eqnarray}\label{e:o2}
\int_{\R_+^n}\d _1^{\frac{n+2}{n-2}} \d _2= \frac{1}{2}c_2\e _{12} + 
O(\e _{12}^{\frac{n}{n-2}}log(\e _{12}^{-1})).
\end{eqnarray}
Thus 
$$
N=\gamma + \a _1\a _2c_2\e _{12} + ||v||^2+ 
O(\e _{12}^{\frac{n}{n-2}}log(\e _{12}^{-1}))
$$
For the denominator, we write
\begin{eqnarray*}
D^{\frac{n}{n-2}}&=&\int_{S_+^n}K (\a _1\tilde\d _1+\a _2
\tilde\d _2)^{\frac{2n}{n-2}}+ 
\frac{2n}{n-2} \int_{S_+^n}K (\a _1\tilde\d _1+\a _2
\tilde\d _2)^{\frac{n+2}{n-2}}v\\
 & + & \frac{n(n+2)}{(n-2)^2}\int_{S_+^n}K(\a _1\tilde\d _1+\a _2
\tilde\d _2)^{\frac{4}{n-2}}v^2\\
 &+& O\left(\int_{S_+^n}(\a _1\tilde\d _1+\a _2
\tilde\d _2)^{\frac{4}{n-2}-1}\inf((\a _1\tilde\d _1+\a _2
\tilde\d _2),|v|)^3+\int |v|^{\frac{2n}{n-2}}\right).
\end{eqnarray*}
We also write
\begin{align*}
\int_{S_+^n}K (\a _1\tilde\d _1+\a _2\tilde\d _2)^{\frac{2n}{n-2}} & = 
\int_{S_+^n}K (\a _1\tilde\d _1)^{\frac{2n}{n-2}}+ 
\int_{S_+^n}K (\a _2\tilde\d _2)^{\frac{2n}{n-2}}+\frac{2n}{n-2}  
\int_{S_+^n}K (\a _1\tilde\d _1)^{\frac{n+2}{n-2}}\a _2\tilde\d _2\\ 
 & + \frac{2n}{n-2}\int_{S_+^n}K (\a _2\tilde\d _2)^{\frac{n+2}{n-2}} 
\a _1\tilde\d _1 +O\left(\int_{S_+^n}\sup(\tilde\d _1,
\tilde\d _2)^{\frac{4}{n-2}}  \inf(\tilde\d _1,\tilde\d _2)^2\right).
\end{align*}
Expansions of $K$ around $a_1$ and $a_2$ give
\begin{eqnarray}\label{e:o3}
\int_{S_+^n}K (\tilde\d _i)^{\frac{2n}{n-2}}=K(a_i)\frac{S_n}{2}-
\frac{2c_1}{\l _i}\frac{\partial K}{\partial \nu}(a_i)+O(\frac{1}{\l
  _i ^2})
\end{eqnarray}
\begin{eqnarray}\label{e:o4}
\int_{S_+^n}K (\tilde\d _i)^{\frac{n+2}{n-2}}\tilde\d _j=K(a_i)
\frac{c_2}{2}\e _{12}+O\left(\e _{12}^{\frac{n}{n-2}}log(\e _{12}^{-1})+
\frac{\e _{12}}{\l _i}(log(\e _{12}^{-1}))^{\frac{n-2}{n}}\right)
\end{eqnarray}
It easy to check
\begin{eqnarray}\label{e:o5}
\int_{S_+^n}\sup^{\frac{4}{n-2}}(\tilde\d _1,\tilde\d _2)\inf^2
(\tilde\d _1,\tilde\d _2)=
O(\e _{12}^{\frac{n}{n-2}}log(\e _{12}^{-1})) \, \mbox{ if } n\geq 4
\end{eqnarray}
and 
\begin{eqnarray}\label{e:o6}
\int_{S_+^n}(\a _1\tilde\d _1+\a _2
\tilde\d _2)^{\frac{4}{n-2}-1}\inf((\a _1\tilde\d _1+\a _2
\tilde\d _2),|v|)^3+\int |v|^{\frac{2n}{n-2}}=O\left(||v||^{\inf(3,
\frac{2n}{n-2})}\right)
\end{eqnarray}
Combining \eqref{e:o1},\eqref{e:o2}, \eqref{e:o3}, \eqref{e:o4},
\eqref{e:o5} and  \eqref{e:o6}, we easily derive our proposition.
\hfill$\Box$\\
A natural improvement of Proposition \ref{p:31} is obtained by taking care 
of the $v$-part, in order to show that it can be neglected with respect to 
the concentration phenomenon.\\
Set 
$$
E_\e =\{v\in H^1(S_+^n) \, /\, ||v||\leq \e \mbox{ and } v 
\mbox{ satisfies } (V_0)\}
$$
where $(V_0)$ is the following condition 
$$
(V_0)\qquad (v,\tilde\d _i)=(v,\partial \tilde\d _i/\partial \l _i)=
(v,\partial \tilde\d _i/\partial a_i)=0, \, \mbox{ for } i=1,2.
$$
Notice that, one can prove arguing as in \cite{B2} (see also \cite{R}), that 
for $\e $ small enough, there exists $\rho >0$ such that for all $v\in E_\e $ 
$$
Q(v,v) \geq \rho ||v||^2.
$$
It follows the following lemma whose proof is similar , up to minor modifications to  corresponding statements in \cite{B2} (see also 
\cite{R}).
\begin{lemma}\label{l:32}
There exists a $C^1$-map which, to each $(\a , a, \l )$ such that 
$ \a _1\tilde\d _1+\a _2\tilde\d _2 \in V_b(2,\e )$ with small $\e $, 
associates $\overline{v}=\overline{v}_{(\a ,a,\l )}$ satisfying 
$$
J( \a _1\tilde\d _1+\a _2\tilde\d _2  +\overline{v})= \min\{
J( \a _1\tilde\d _1+\a _2\tilde\d _2  +v) , \, v \mbox{ satisfies } (V_0)\}.
$$
Moreover, there exists $c>0 $ such that the following holds
$$
||\overline{v}||\leq c \left(\frac{1}{\l _1}+\frac{1}{\l _2}+\e _{12}^{\frac 
{n+2}{2(n-2)}}log(\e _{12}^{-1})+(\mbox{ if } n\leq 5) \e _{12}
(log(\e _{12}^{-1}))^{\frac{n-2}{n}}\right).
$$
\end{lemma}
\begin{proposition}\label{p:33}
Let $n\geq 4$, for $ u=\a _1\tilde\d _1+\a _2\tilde\d _2 \in V_b(2,\e )$, we 
have the following expansion 
\begin{align*}
(\n J(u),\l _1\partial\tilde\d _1/\partial\l _1)= & 2J(u)\left[\frac{c_2}
{2}\a _2 \l _1\frac{\partial \e _{12}}{\partial \l _1}(1-J(u)^\frac{n}{n-2}
(\a _1^\frac{4}{n-2}K(a_1)+\a _2^\frac{4}{n-2}K(a_2)))\right.\\
 & \left.-2J(u)^\frac{n}{n-2}\a _1^\frac{n+2}{n-2}\frac{c_3}{\l _1}\frac{\partial
K}{\partial \nu}(a_1)\right]+O(\frac{1}{\l _1^2})\\
 & +O\left(\e _{12}^{\frac{n}{n-2}}log(\e _{12}^{-1})+\e _{12}(log(
\e _{12}^{-1}))^{\frac{n-2}{n}}(\frac{1}{\l _1}+\frac{1}{\l _2})\right).
\end{align*}
where $c_3=\frac{n-2}{2}c_0^{2n/(n-2)}\int_{\R_+^n}\frac{x_n(|x|^2-1)}
{(1+|x|^2)^{n+1}}dx $
\end{proposition} 
\proof
We have 
\begin{eqnarray}\label{e:31}
(\n J(u),h)=  2J(u)\left[\int_{S_+^n}\n u\n h+\frac{n(n-2)}{4}
\int_{S_+^n}uh-J(u)^\frac{n}{n-2}\int_{S_+^n}Ku^\frac{n+2}{n-2}h\right]
\end{eqnarray}
Observe that(see \cite{B2})
\begin{eqnarray}\label{e:32}
\int _{\R^n_+}\n \d _1 \n (\l _1 \frac{\partial \d _1}{\partial \l
  _1})
= \int _{\R^n_+}\d _1^{\frac{n+2}{n-2}}\l _1\frac{\partial \d _1}{\partial \l
  _1} = 0
\end{eqnarray}
\begin{align}\label{e:33}
 \int _{\R^n_+}\n \d _2 \n (\l _1 \frac{\partial \d _1}{\partial \l
  _1})
&= \int _{\R^n_+}\d _2^{\frac{n+2}{n-2}}\l _1\frac{\partial \d _1}{\partial \l
  _1}\nonumber \\
 &=\frac{1}{2} c_2 \l _1 \frac{\partial \e _{12}}{\partial \l _1} +
  O\left( \e _{12}^{\frac{n}{n-2}}log(\e _{12}^{-1})\right) 
\end{align}
\begin{align}\label{e:34}
\int _{\R^n_+}K \d _1 ^{\frac{n+2}{n-2}} \l _1 \frac{\partial \d
  _1}{\partial \l _1}
&= 2 \n K(a_1)\int _{\R^n_+}\d _1 ^{\frac{n+2}{n-2}}\l _1 \frac{\partial
  \d _1}{\partial \l _1} (x-a_1) +O\left(\frac{1}{\l _1
  ^2}\right)\nonumber\\
&=- \frac{2c_3}{\l _1}\n K(a) e_n + O\left(\frac{1}{\l _1 ^2}\right),
\end{align}
\begin{align}\label{e:35}
\int _{\R^n_+}K\d _2 ^{\frac{n+2}{n-2}}\l _1 \frac{\partial \d
  _1}{\partial \l _1}& = K(a_2) \frac{1}{2}c_2 \l _1 \frac{\partial
  \e _{12}}{\partial \l _1} +
  O\left(\frac{1}{\l _2}\e _{12} (Log(\e _{12}
  ^{-1}))^{\frac{n-2}{n}}\right)\nonumber\\
& + O\left( \e _{12}^{\frac{n}{n-2}}log(\e _{12}^{-1})\right),
\end{align}
\begin{align}\label{e:36}
\frac{n+2}{n-2}\int _{\R^n_+}K\d _2\d _1 ^{\frac{4}{n-2}}\l _1 \frac{\partial \d _1}{\partial \l _1}& = K(a_1) \frac{1}{2}c_2 \l _1 \frac{\partial
  \e _{12}}{\partial \l _1} +
  O\left(\frac{1}{\l _1}\e _{12} (Log(\e _{12}
  ^{-1}))^{\frac{n-2}{n}}\right)\nonumber\\
& + O\left(\e _{12}^{\frac{n}{n-2}}log(\e _{12}^{-1})\right).
\end{align}
Combining \eqref{e:31}, \eqref{e:32}, \eqref{e:33}, \eqref{e:34},
\eqref{e:35} and \eqref{e:36}, we easily derive our proposition.
\hfill$\Box$\\
\begin{proposition}\label{p:34}
Let $n\geq 4$. For $ u = \sum \a _i \tilde\d _i \in V_b (2, \e )$, we
have the following expansion:
\begin{align*}
\left(\n J(u), \frac{1}{\l _1}\frac{\partial\tilde\d _1}{\partial a_1}\right)
&= 2J(u)\a _1 e_n \left[ c_4\left(1-
      J(u)^{\frac{n}{n-2}}\a _1 ^{\frac{4}{n-2}}K(a_1)\right)+
    J(u)^{\frac{n}{n-2}}\a _1 ^{\frac{4}{n-2}}\frac{c_5}{\l
      _1}\frac{\partial K}{\partial\nu}(a_1)\right]\\
&-J(u)\a _2 c_2 \frac{1}{\l
  _1}\frac{\partial\e _{12}}{\partial a_1}\left(-1 +
  J(u)^{\frac{n}{n-2}}\sum \a _i ^{\frac{4}{n-2}}K(a_i)\right)\\
&-4J(u)^{\frac{2(n-1)}{n-2}}\a _1 ^{\frac{n+2}{n-2}}\frac{2c_5}{\l _1}\n
_T K(a_1) + O\left(\e _{12}^{\frac{n}{n-2}}Log(\e _{12}^{-1})+ \e
  _{12}^{\frac{n+1}{n-2}}\l _2 |a_1-a_2|\right)\\
&+O\left(\e _{12}\left(Log(\e _{12}^{-1})\right)^{\frac{n-2}{n}}\sum
  \frac{1}{\l _k}\right)+ O\left(\frac{1}{\l _1^2}\right)
\end{align*}
where
$$
c_4 = (n-2)c_0^{\frac{2n}{n-2}}\int
_{\R^n_+}\frac{x_n}{(1+|x|^2)^{n+1}}dx \quad \mbox{and }
c_5=\frac{n-2}{2n}c_0^{\frac{2n}{n-2}}\int
_{\R^n}\frac{x_n^2}{(1+|x|^2)^{n+1}}dx
$$
\end{proposition}
\proof
An easy computation shows
\begin{eqnarray}\label{e:37}
\int _{\R^n_+}\n\d _1\n\left(\frac{1}{\l _1}\frac{\partial\d
      _1}{\partial a_1}\right)=\int _{\R^n_+}\d
      _1^{\frac{n+2}{n-2}}\frac{1}{\l _1}\frac{\partial\d _1}{\partial
      a_1} = c_4e_n,
\end{eqnarray}
\begin{align}\label{e:38} 
  \int _{\R^n_+}\n\d _2\n\left(\frac{1}{\l _1}\frac{\partial\d
      _1}{\partial a_1}\right) 
&=\int _{\R^n_+}\d
      _2^{\frac{n+2}{n-2}}\frac{1}{\l _1}\frac{\partial\d _1}{\partial
      a_1}\nonumber\\
&= \frac{1}{2}\frac{c_2}{\l _1}\frac{\partial\e _{12}}{\partial a_1} + O\left(\e
      _{12}^{\frac{n}{n-2}}Log(\e _{12}^{-1}) +
       \e _{12}^{\frac{n+1}{n-2}}\l _2 |a_1-a_2|\right)
\end{align}
\begin{eqnarray}\label{e:39}
\int _{\R^n_+}K\d
      _1^{\frac{n+2}{n-2}}\frac{1}{\l _1}\frac{\partial\d _1}{\partial
      a_1} =K(a_1) c_4e_n + 2\frac{c_5}{\l _1}\n K(a_1)
      +O\left(\frac{1}{\l _1^2}\right),
\end{eqnarray}
\begin{align}\label{e:310}
\int _{\R^n_+}K\d
      _2^{\frac{n+2}{n-2}}\frac{1}{\l _1}\frac{\partial\d _1}{\partial
      a_1} &=K(a_2)\frac{1}{2}c_2\frac{1}{\l _1}\frac{\partial\e
      _{12}}{\partial a_1}+ O\left(\e _{12}^{\frac{n+1}{n-2}}\l
      _2|a_1-a_2|\right)\nonumber\\
  &+O\left(\e _{12}^{\frac{n}{n-2}}Log(\e
      _{12}^{-1})\right)+O\left(\frac{1}{\l _2}\e _{12}\left(Log(\e
      _{12}^{-1})\right)^{\frac{n-2}{n}}\right),
\end{align}
\begin{align}\label{e:311}
\frac{n+2}{n-2}\int _{\R^n_+}K\d
      _1^{\frac{4}{n-2}}\d _2 \frac{1}{\l _1}\frac{\partial\d _1}{\partial
      a_1} &=K(a_1)\frac{1}{2}c_2\frac{1}{\l _1}\frac{\partial\e
      _{12}}{\partial a_1}+ O\left(\e _{12}^{\frac{n+1}{n-2}}\l
      _2|a_1-a_2|\right),\nonumber\\
  &+O\left(\e _{12}^{\frac{n}{n-2}}Log(\e
      _{12}^{-1})\right)+O\left(\frac{1}{\l _1}\e _{12}\left(Log(\e
      _{12}^{-1})\right)^{\frac{n-2}{n}}\right).
\end{align}
Using \eqref{e:31}, \eqref{e:37}, \eqref{e:38}, \eqref{e:39},
\eqref{e:310} and \eqref{e:311}, our proposition follows.
\hfill$\Box$\\
\begin{proposition}\label{p:35}
 For $ u = \sum \a _i \tilde\d _i \in V_b (2, \e )$, we
have the following expansion:
\begin{align*}
\left(\n J(u),\tilde\d _1\right)&=2J(u) \a _1
  \frac{S_n}{2}\left(1-J(u)^{\frac{n}{n-2}}\a _1^{\frac{4}{n-2}}K(a_1)\right)
+2J(u)^{\frac{2(n-1)}{n-2}}\a _1^{\frac{n+2}{n-2}}\frac{2c_6}{\l _1}
\frac{\partial K}{\partial \nu}(a_1)\\
&-J(u)c_2\e _{12}\a _2 \left(-1 +
  J(u)^{\frac{n}{n-2}}\left(\frac{n+2}{n-2}\a
  _1^{\frac{4}{n-2}}K(a_1)+\a _2^{\frac{4}{n-2}}K(a_2)\right)\right)\\
&+O\left(\frac{1}{\l _1^2}+ \e _{12}^{\frac{n}{n-2}}Log(\e _{12}^{-1})
  + \frac{1}{\l _2}\e _{12}\left(Log(\e _{12}^{-1})\right)^{\frac{n-2}{n}}\right)
\end{align*}
where $c_6= c_0 \int _{\R^n_+}\frac{x_n}{(1+|x|^2)^n}dx$ and where
$S_n$ is defined in Proposition \ref{p:31}.
\end{proposition}
\proof
Using estimates \eqref{e:o1}, \eqref{e:o2}, \eqref{e:o3} and
\eqref{e:o4}, we easily derive our proposition.
\hfill$\Box$\\
Before ending this section, we state  the above results in the case
where we  have only one mass instead of two masses.

\begin{proposition}\label{p:36}
For $\e > 0$ small enough and $ u=\a\tilde\d _{a,\l } + v \in V_b(1,\e
) $, we have the following expansion:
\begin{align*}
J(u)=&
\frac{(S_n/2)^{(2/n)}}{(K(a))^{\frac{n-2}{n}}}\left[1+\frac{8c_1(n-2)}
{nK(a)S_n \l}
\frac{\partial K}{\partial 
\nu}(a) +f(v)\right.\\
 &\left. +Q(v,v) +O\left(\frac{1}{\l ^2}\right)+ O\left(
||v||^{\inf(3,\frac{2n}{n-2})}\right)\right]
\end{align*} 
where 
\begin{align*}
Q(v,v) = & \frac{2}{\a ^2 S_n}\left[||v||^2-\frac{n+2}{(n-2) K(a)} \int_{S_+^n}K(\tilde\d _{a,\l })^{\frac{4}{n-2}}v^2\right]  \\ 
f(v) = & -\frac{4}{\a K(a)S_n }\int_{S_+^n}K(\tilde\d _{a,\l })^{\frac{n+2}{n-2}}v.
\end{align*}
\end{proposition} 
\begin{proposition}\label{p:37}
For $u=\a \tilde\d _{a,\l} \in V_b(1,\e)$, we have the following
expansion:
\begin{eqnarray*}
\left(\n J(u), \l
  \frac{\partial\tilde\d}{\partial\l}\right)=-4J(u)^{\frac{2(n-1)}{n-2}}\a
  ^{\frac{n+2}{n-2}}\frac{c_3}{\l}\frac{\partial K}{\partial\nu}(a) +
  O\left(\frac{1}{\l ^2}\right)
\end{eqnarray*}
\end{proposition}
\begin{proposition}\label{p:38}
For $u=u=\a \tilde\d _{a,\l} \in V_b(1,\e)$, we have the following
expansion:
\begin{align*}
\left(\n J(u),\frac{1}{ \l}
  \frac{\partial\tilde\d}{\partial a}\right)&=2J(u)\a
  e_n\left[c_4\left(1-J(u)^{\frac{n}{n-2}}\a
  ^{\frac{4}{n-2}}K(a)\right)+2c_5J(u)^{\frac{n}{n-2}}\frac{\a
  ^{\frac{4}{n-2}}}{\l}\frac{\partial K}{\partial\nu}(a)\right]\\
&-4J(u)^{\frac{2(n-1)}{n-2}}\a ^{\frac{n+2}{n-2}}c_5 \frac{\n _T K(a)}{\l}
  + O\left(\frac{1}{\l ^2}\right).
\end{align*}
\end{proposition}

\section{Morse Lemma at Infinity}

\mbox{}
In this section, we consider the case where we only have one mass and
we perform a Morse lemma at infinity for J, which completely gets rid
of the v-contribution and shows that the functional behaves, at
infinity, as $J(\a \tilde\d _{\tilde{a},\tilde{\l}}) + |V|^2$, where $V$ is a variable
completely independent of $\tilde{a}$, $\tilde{\l}$. Namely, we prove the following
proposition.
\begin{proposition}\label{p:41}
For $\e > 0$ small enough, there is a diffeomorphism
$$
\a \tilde\d _{a,\l} + v \longmapsto \a \tilde\d _{\tilde{a},\tilde\l
  }\in V_b(1,\e ' )
$$
for some $\e '$ such that
$$
J(\a \tilde\d _{a,\l}+v)= J(\a \tilde\d _{\tilde{a}, \tilde\l})+ |V|^2
$$
where $V$ belongs to a neighborhood of zero in a fixed Hilbert space.
\end{proposition}
The above Morse Lemma can be improved when the concentration point is
near a critical point $y$ of $K_1=K_{/\partial S^n_+}$ with
$\frac{\partial K}{\partial \nu}(y) > 0 $ , leading to the following
normal form:
\begin{proposition}\label{p:42}
For $u=\a \tilde\d _{\tilde{a},\tilde\l}\in V_b(1,\e )$ such that
$\tilde{a} \in \mathcal{V}(y,\rho )$, $ \frac{\partial
  K}{\partial\nu}(y) > 0$, $\rho > 0$ and $y$ is a critical point of
$K_1$, there is another change of variable $(\tilde{\tilde{a}},
\tilde{\tilde\l} )$ such that
\begin{eqnarray*}
J(u)&=\psi (\tilde{\tilde{a}}, \tilde{\tilde{\l}})\\
&:= \frac{(S_n/2)^{\frac{2}{n}}}{(K(\tilde{\tilde{a}}))^{\frac{n-2}{n}}}\left[1+(c-\eta
  )\frac{1}{\tilde{{\tilde{\l}}}} \frac{\partial K}{\partial\nu} (\tilde{\tilde{a}}) \, \right]
\end{eqnarray*}
where $\eta$ is a small positive real.\\
Here and in the sequel, $\mathcal{V}(y,\rho )$ denotes a neighborhood
of $y$.
\end{proposition}
The proof of Propositions \ref{p:41} and \ref{p:42} can be easily
deduced from the following
lemma, arguing as in \cite{B1} and \cite{BCCH} .
\begin{lemma}\label{l:43}
There exists a pseudogradient $Z$ so that the following holds.\\
There is a constant $c>0$ independent of $u=\a\tilde\d _{a,\l }$ in
$V_b(1,\e )$ such that \\
{\bf i.}\quad $-\left(\n J(u), Z\right) \geq \frac{c}{\l}$\\
{\bf ii.}\quad $-\left(\n J(u+\bar{v}), Z +\frac{\partial
    \bar{v}}{\partial (\a , a, \l )}(Z)\right) \geq \frac{c}{\l}$\\
{\bf iii.} \quad $ Z$ is bounded\\
{\bf iv.}\quad the only region where $\l$ increases along $Z$ is the
region where $a\in \mathcal{V}(y,\rho )$, where $y$ is a critical
point of $K_1$ such that $\frac{\partial K}{\partial\nu}(y) > 0. $
\end{lemma}
Before giving the proof of Lemma \ref{l:43}, we notice that combining
Proposition \ref{p:42} and Lemma \ref{l:43}, one can easily derive the
following corollary.
\begin{corollary}\label{c:44}
Assume that $J$ does not have any critical point. Then, the only
critical points at infinity of $J$ in $V_b(1,\e )$, for $\e $ small
enough, correspond to $\tilde\d _{y,\infty}$, where $y$ is a critical
point of $K_1 = K_{/\partial S^n_+}$ such that $\frac{\partial
  K}{\partial\nu}(y) > 0 $.\\
Moreover such a critical point at infinity has a Morse index equal to
$(n-1-index(K_1,y))$.
\end{corollary}
\begin{pfn}{\bf Lemma \ref{l:43}}
Let $u=\a \tilde\d _{a,\l} \in V_b(1,\e )$. We divide $V_b(1,\e )$ in
three regions.\\
{\bf $1^{st}$ region}. \quad $a\notin \cup _{0\leq i\leq s} \mathcal{V}(y_i,
\rho )$, where $\rho < \frac{1}{2} \min _{i\not= j}d(y_i,y_j)$.\\
Set $Z_1=\frac{1}{\l }\frac{\partial\tilde\d }{\partial a}\n _T K(a)$,
from Proposition \ref{p:38}, we have
\begin{align*}
-(\n J(u), Z_1)&=c\frac{|\n _TK(a)|^2}{\l} + O\left(\frac{1}{\l ^2}\right)\\
               &\geq \frac{c}{\l}
\end{align*}
{\bf $2^{nd}$ region}. \quad $a\in \cup _{0\leq i\leq l} \mathcal{V}(y_i,
2\rho )$\\
We set 
$$
Z_2 = \l \frac{\partial\tilde\d}{\partial\l} +
\frac{1}{\l}\frac{\partial\tilde\d}{\partial a} \n _TK(a).
$$
Using Propositions \ref{p:37} and \ref{p:38}, we obtain\\
\begin{align*}
-(\n J(u),Z_2) &\geq c \frac{|\n
  _TK(a)|^2}{\l}+\frac{c}{\l}\frac{\partial K}{\partial\nu}(a)+
O\left(\frac{1}{\l ^2}\right)\\
&\geq \frac{c}{\l}.
\end{align*}
{\bf $3^{rd}$ region}. \quad $a\in \cup _{l+1\leq i\leq s} \mathcal{V}(y_i,
2\rho ).$\\
We set 
$$
Z_3 = -\l \frac{\partial\tilde\d}{\partial\l}. 
$$
Using Proposition \ref{p:37}, we deduce that
\begin{align*}
-(\n J(u), Z_3) &\geq -\frac{c}{\l}\frac{\partial K}{\partial\nu}(a) +
 O\left(\frac{1}{\l ^2}\right)\\
&\geq \frac{c}{\l}.
\end{align*}
Hence our global vector field will be built using a convex combination
of $Z_1$, $Z_2$ and $Z_3$ and will satisfy obviously i., iii. and
iv. Regarding the estimate ii., it can be obtained once we have
i. arguing as in \cite{B1} and \cite{BCCH}.
\end{pfn}

\section{Proof of Theorem \ref{t:11}}

\mbox{}
Our proof follows the algebraic topological arguments introduced in \cite{AB}. Arguing by contradiction, we suppose that $J$ has no critical points. It follows
from Corollary \ref{c:44}, that under the assumptions of Theorem
\ref{t:11}, the critical points at infinity of $J$ under the level
$c_1= (S_n/2)^{\frac{2}{n}}(K(y_l))^{\frac{2-n}{n}} + \e $ 
, for $\e$ small enough, are in one to one correspondance with the  critical points of $K_1$ $y_0$, $y_1$,
..., $y_l$. The unstable manifold at infinity of such  critical
points at infinity, $W_u(y_0)_\infty$, ..., $W_u(y_l)_\infty$ can be
described, using Proposition \ref{p:42}, as the product of $W_s(y_0)$,
..., $W_s(y_l)$ (for a pseudogradient of $K$ ) by $[A, +\infty [$
 domaine of the variable $\l$, for some positive
number $A$ large enough.\\
Since $J$ has no critical points, it follows that $ J_{c_1}=\{u\in
\sum ^+ / J(u) \leq c_1 \}$ retracts by deformation on $X_\infty =
\cup _{0\leq j\leq l}W_u(y_j)_\infty$ (see Sections 7 and 8 of
\cite{BR}) which can be parametrized as we said before by $X \times
[A, +\infty[$.\\
From another part, we have $X_\infty$ is contractible in $J_{c_2+\e}$,
where $c_2=(S_n/2)^{\frac{2}{n}}c^{\frac{2-n}{n}}$. Indeed from
$(A_4)$, it follows that there exists a contraction $ h :[0,1] \times
X \to K^c$, $h$ continuous such that for any $a\in X \quad h(0,a)=a$
and $h(1,a)=a_0$ a point of $X$. Such a contraction  gives rise to
the following contraction $\tilde{h} : X_\infty \to \sum ^+$ defined
by 
$$
[0,1] \times X \times \left[0,\right.+\infty\left[ \right.\ni (t,a_1,\l _1 ) \longmapsto \tilde\d _{(h(t,a_1),\l )}
+ \bar{v} \in \Sigma ^+, \quad a_1 \in X, \, \, \l _1 \geq A
$$
For $t=0$, $\tilde\d _{(h(0,a_1),\l _1)}+\bar{v} = \tilde\d _{a_1, \l _1}
+\bar{v} \in X_\infty$. $\tilde{h}$ is continuous and
$\tilde{h}(1,a_1,\l _1)= \tilde\d _{a_0,\l _1} +\bar{v}$, hence our
claim follows.\\
Now, using Proposition \ref{p:36}, we deduce that
$$
J(\tilde\d _{h(t,a_1), \l _1} + \bar{v}) \sim
(\frac{S_n}{2})^{\frac{2}{n}}(K(h(t,a_1)))^{\frac{2-n}{n}}\left(1+O(A
    ^{-2})\right)
$$
where $K(h(t,a_1)) \geq c $ by construction.\\
Therefore such a contraction is performed under $c_2 +\e$, for $A$
large enough, so $X_\infty$ is contractible in $J_{c_2+\e }$.\\
In addition, choosing $c_0$ small enough, $J_{c_2+\e }$ retracts by
deformation on $J_{c_1}$, which retracts by deformation on $X_\infty$,
therfore $X_\infty$ is contractible leading to the contractibility of
$X$, which is in contradiction with our assumption. Hence our theorem
follows.\\\\
Before ending this section, we give the proof of Corollary
\ref{c:12}.\\ 

\begin{pfn}{\bf Corollary \ref{c:12}}
Arguing by contradiction, we may assume that the Morse index of the
solution provided by Theorem \ref{t:11} is $\leq m-1$.\\
Perturbing, if necessary $J$, we may assume that all the critical
points of $J$ are nondegenerate and have their Morse index $\leq
m-1$. Such critical points do not change the homological group in
dimension $m$ of level sets of $J$.\\
Since $X_\infty$ defines a homological class in dimension $m$ which is
nontrivial in $J_{c_1}$, but trivial in $J_{c_2+\e}$, our result
follows.
\end{pfn}

\section{Proof of Theorems \ref{t:13} and \ref{t:14}}

\mbox{}
First, we start by proving the following  main results
\begin{proposition}\label{p:61}
Let $y_0$ be defined in $(H_1)$. Then $(y_0,y_0)_\infty$
is not a critical point at infinity for $J$, that is, there exists a
decreasing pseudogradient $W$ for $J$ satisfying Palais-Smaile in the
neighborhood of $(y_0,y_0)_\infty$. 
\end{proposition}
\proof
For $\e _0 > 0$ small enough, we set
$$
C_{\e _0}=\{u=\a _1 \tilde\d _{a_1, \l _1} + \a _2 \tilde\d _{a_2,\l
  _2} \in V_b (2,\e _0 )/ a_1, a_2 \in \mathcal{V}(y_0)\cap \partial
S^n_+ \}.
$$
Our goal is to build a pseudogradient vector field $W$ for $J$ satisfying the 
Palais-Smale condition in $C_{\e _0}$ such that for $u\in C_{\e _0}$,
we have\\ 
{\bf i.}\quad $-\left(\n J(u), W\right) \geq \gamma \left(\sum\frac{1}{\l
    _i^2}+ \sum (1-J(u)^{\frac{n}{n-2}}a
  _i^{\frac{4}{n-2}}K(a_i)\right) + c\e _{12}^{\frac{n-\frac{1}{2}}{n-2}},$\\
{\bf ii.}\quad $-\left(\n J(u+\bar{v}), W +\frac{\partial
    \bar{v}}{\partial (\a , a, \l )}(W)\right) \geq
\frac{\gamma}{2} \left(\sum\frac{1}{\l
    _i^2}+ \sum (1-J(u)^{\frac{n}{n-2}}a
  _i^{\frac{4}{n-2}}K(a_i)\right) + c\e _{12}^{\frac{n-\frac{1}{2}}{n-2}}.$\\
{\bf iii.} \quad $ W$ is bounded\\
{\bf iv.} \quad $\dot{\l}_{max} \leq 0$.\\
where $\gamma$  is a positive constant large enough.\\
We can assume, without loss of generality, that $\l _1 \leq \l _2$. We
devide $C_{\e _0}$ in three principal regions.\\
{\bf $1^{st}$ region.}\quad $M\l _1 \leq \l _2$ and $\forall i |1- J(u)^{\frac{n}{n-2}}a
  _i^{\frac{4}{n-2}}K(a_i)| \leq \frac{2C'}{\l _i}$, where $M$ and
  $C'$ are 
postive constants large enough.\\
We set
$$
W_1 = \l _1 \frac{\partial\tilde\d _1}{\partial\l _1}- \l _2
\frac{\partial\tilde\d _2}{\partial\l _2}\sqrt{M}.
$$
From Proposition \ref{p:33}, we derive
\begin{align*}
-\left(\n J(u), W_1\right)&\geq c\left(\frac{1}{\l _1} + O(\e
 _{12})\right) + \sqrt{M}\left[-\frac{c}{\l _2} + c \e _{12}\right]
 +O\left(\sum \frac{1}{\l _k^2}\right)\\
&+O\left(\e _{12}^{\frac{n}{n-2}}Log(\e _{12}^{-1}) + \e _{12}\left(Log(\e
 _{12}^{-1})\right)^{\frac{n-2}{n}}\sum \frac{1}{\l _k}\right)\\
&\geq \frac{c}{2}\left(\frac{1}{\l _1}+ \frac{M}{\l _2}\right) +
 \frac{c\sqrt{M}}{2}\e _{12}\\
&\geq C\left(\sum\frac{1}{\l
    _i^2}+ \sum (1-J(u)^{\frac{n}{n-2}}a
  _i^{\frac{4}{n-2}}K(a_i) + \e _{12}\right).
\end{align*}
{\bf $2^{nd}$ region.}\quad $2M\l _1 \geq \l _2$ and, $\forall i |1- J(u)^{\frac{n}{n-2}}a
  _i^{\frac{4}{n-2}}K(a_i)| \leq \frac{2C'}{\l _i}$\\
In this region, two cases may occur.\\
{\bf Case 2.1}\quad $\e _{12}^{\frac{n-\frac{1}{2}}{n-2}} \geq
\frac{d_0}{\l _1}$, where $d_0=\max (d(y_0,a_1), d(y_0,a_2),
d(a_1,a_2)).$\\
We set
$$
W_2=\frac{1}{\l _1}\left[\a _1 \frac{\partial\tilde\d _1}{\partial
    a_1}\left(\frac{a_2-a_1}{d(a_2,a_1)}\right) -\a _2
  \frac{\partial\tilde\d _2}{\partial
    a_2}\left(\frac{a_2-a_1}{d(a_2,a_1)}\right)\right].
$$
From Proposition \ref{p:34} and the fact that $\frac{\partial \e
  _{12}}{\partial a_2} = -\frac{\partial \e
  _{12}}{\partial a_1} $, we obtain
 \begin{align*}
-\left(\n J(u), W_2\right)&=\frac{1}{\l _1}\left(\a _1\a _2 c_2
 2J(u) \frac{\partial\e _{12}}{\partial
 a_1}\right)\left(-1 +J(u)^{\frac{n}{n-2}}\sum \a
 _k^{\frac{4}{n-2}}K(a_k)\right)2\left(\frac{a_2-a_1}{d(a_2,a_1)}\right)\\
&+O\left(\sum \left(\frac{|\n _TK(a_k)|}{\l _k} + \frac{1}{\l _k^2}\right)\right)\\
&+O\left(\e _{12}^{\frac{n}{n-2}}Log(\e _{12}^{-1}) + \sum \frac{1}{\l
 _k} \e
 _{12}\left(Log(\e _{12}^{-1})\right)^{\frac{n-2}{n}}+ 
 d(a_1,a_2)\e _{12}^{\frac{n+1}{n-2}} \sum \l _k \right).
\end{align*}
Observe that
\begin{align*}
\e _{12} &\sim \left(\l _1\l _2 d(a_1,a_2)^2 \right)^{\frac{2-n}{2}}
\quad \mbox.{Indeed } \l _1 \mbox{and } \l _2 \mbox{ are the same order}, and\\
\frac{\partial \e _{12}}{\partial a_1}&= -(n-2)\l _1\l _2 (a_1-a_2)\e
_{12}^{\frac{n}{n-2}}.
\end{align*}
Thus
\begin{align}\label{e:41}
 -\left(\n J(u), W_2\right)&=\frac{4}{\l _1}\a _1\a _2 c_2 (n-2)
 J(u)\l _1\l _2 d(a_1,a_2)\e
 _{12}^{\frac{n}{n-2}}(1+o(1)) + R\nonumber\\
&\geq \frac{c\e_{12}}{\l _1 d(a_1,a_2)} +R\nonumber\\
&\geq c\e _{12}^{\frac{n-1}{n-2}} +R
\end{align}
where 
\begin{align*}
R&= +O\left(\frac{d(a_1,y_0)}{\l _1}+ \frac{d(a_2,y_0)}{\l _2} +\sum
 \frac{1}{\l _k^2} \right)\\
&+O\left(\e _{12}^{\frac{n}{n-2}}Log(\e _{12}^{-1}) + \sum \frac{1}{\l
 _k} \e
 _{12}\left(Log(\e _{12}^{-1})\right)^{\frac{n-2}{n}}+ 
 d(a_1,a_2)\e _{12}^{\frac{n+1}{n-2}} \sum \l _k\right).
\end{align*}
We  also observe  that
\begin{eqnarray}\label{e:42}
\frac{d_0}{\l _1} \leq \e _{12}^{\frac{n-\frac{1}{2}}{n-2}}= o\left(\e
  _{12}^{\frac{n-1}{n-2}}\right),
\end{eqnarray}
\begin{eqnarray}\label{e:43}
\frac{1}{\l _1d(a_1,a_2)}\geq \left(\frac{1}{\l _1\l _2
    d(a_1,a_2)^2}\right)^{\frac{1}{2}} \sim \e _{12}^{\frac{1}{n-2}},
\end{eqnarray}
\begin{eqnarray}\label{e:44}
\l _id(a_1,a_2)\e _{12}^{(n+1)/(n-2)}= \e _{12}^{n/(n-2)}\left(\frac{\l
  _i}{\l _j}\right)^{1/2}=o(\e _{12}^{(n-1)/(n-2)}),
\end{eqnarray}
\begin{align}\label{e:45}
\frac{\e _{12}}{\l _1}\left(Log(\e _{12}^{-1})\right)^{(n-2)/n}&=
\frac{\e _{12}}{\sqrt{\l _1d_0}}\frac{\sqrt{d_0}}{\sqrt{\l
    _1}}\left(Log(\e _{12}^{-1})\right)^{(n-2)/n}\nonumber\\
&=O\left(\frac{d_0}{\l _1}+ \frac{\e _{12}^2\left(Log(\e
      _{12}^{-1})\right)^{2(n-2)/n} }{\l _1 d_0}\right)\nonumber\\
&=o\left(\e _{12}^{(n-1)/(n-2)}\right)+o\left(\frac{\e _{12}}{\l _1
    d_0}\right).
\end{align}
In the same way, we have
 \begin{eqnarray}\label{e:46}
\frac{\e _{12}}{\l _2}\left(Log(\e _{12}^{-1})\right)^{(n-2)/n}=
o\left(\e _{12}^{(n-1)/(n-2)}\right)+o\left(\frac{\e _{12}}{\l _2
    d_0}\right)
\end{eqnarray}
We also have, since $\l _1 |a_1-a_2| \to +\infty $,
\begin{eqnarray}\label{e:47}
\l _1^{-2}= o\left(d_0 \l _1^{-1}\right)= o\left(\e
  _{12}^{(n-1)/(n-2)}\right).
\end{eqnarray}
Similary, we have
\begin{eqnarray}\label{e:48}
\l _2^{-2}= o\left(\e
  _{12}^{(n-1)/(n-2)}\right).
\end{eqnarray}
Using \eqref{e:41}, \eqref{e:42}, \eqref{e:43}, \eqref{e:44},
\eqref{e:45}, \eqref{e:46}, \eqref{e:47} and  \eqref{e:48}, we find
\begin{align*}
-\left(\n J(u), W_2\right) &\geq C \left(\e _{12}^{(n-1)/(n-2)}
  +\frac{d_0}{\l _1}+ \frac{d_o}{\l _2}\right)\\
&\gamma \left(\e _{12}^{\frac{n-\frac{1}{2}}{n-2}}+\sum \frac{1}{\l
  _i^2}+ \sum |1-J(u)^{\frac{n}{n-2}}\a _i^{\frac{4}{n-2}}K(a_i)^2|\right)
\end{align*}
where $\gamma$ is a large constant.\\
{\bf Case 2.2}\quad $\e _{12}^{\frac{n-\frac{1}{2}}{n-2}} \leq
2\frac{d_0}{\l _1}$\\
In this case, we set
$$
W_3=\frac{1}{\l _1}\left[\a _1 \frac{\partial\tilde\d _1}{\partial
    a_1}\left(\frac{y_0-a_1}{d_o}\right) + \a _2
  \frac{\partial\tilde\d _2}{\partial
    a_2}\left(\frac{y_0-a_2}{d_0}\right)\right].
$$
Using Proposition \ref{p:34}, we obtain
\begin{align*}
-\left(\n J(u), W_3\right)&= J(u)2c_2 \frac{\a _1\a
 _2}{\l _1}\frac{\partial\e _{12}}{\partial
 a_1}\left(\frac{a_2-a_1}{d_0}\right)(1+o(1))\\
&+J(u)^{\frac{2(n-1)}{n-2}}\frac{c_5}{\l _1}\sum \a
 _k^{\frac{2n}{n-2}}\n _TK(a_k)\frac{y_0-a_k}{d_0} +R_1
\end{align*}
where
\begin{align*}
R_1=& O\left(
 \frac{1}{\l _1^2}+ \frac{1}{\l _2^2} + \e _{12}^{\frac{n}{n-2}}Log(\e
 _{12}^{-1}) \right)\\
&+O\left(\sum \frac{1}{\l _k}  \e _{12}\left(Log(\e _{12}^{-1})\right)^{\frac{n-2}{n}}+ 
 d(a_1,a_2)\e _{12}^{\frac{n+1}{n-2}} \sum \l _k\right).
\end{align*}
Thus, we find
\begin{align}\label{e:49}
-\left(\n J(u), W_3\right) & \geq \frac{c}{\l _1}\l _1\l _2
 \frac{d(a_1,a_2)^2}{d_0}\e _{12}^{\frac{n}{n-2}}+\frac{c}{\l
 _1d_0}\left(d(a_1,y_0)^2+d(a_2,y_0)^2\right)+R_1\nonumber\\
&\geq \frac{c}{\l _1d_0}\e _{12}+\frac{c}{\l _1}d_0 + R_1.
\end{align}
Observe that
\begin{eqnarray}\label{e:410}
 \e _{12}^{\frac{n}{n-2}}Log(\e _{12}^{-1}) = o\left(\e
 _{12}^{\frac{n-\frac{1}{2}}{n-2}}\right)= o\left(d_0 \l
 _1^{-1}\right)
\end{eqnarray}
\begin{eqnarray}\label{e:411}
\l _k^{-2}=o\left(d_0 \l _k^{-1}\right)\qquad \mbox{for
  }k=1, 2,
\end{eqnarray}
\begin{align}\label{e:412}
\e _{12}\left(Log(\e _{12}^{-1})\right)^{\frac{n-2}{n}}\frac{1}{\l
  _1}&= \frac{\e _{12}\left(Log(\e
  _{12}^{-1})\right)^{\frac{n-2}{n}}}{\sqrt{\l _1d_0}}
  \frac{\sqrt{d_0}}{\sqrt{\l _1}}\nonumber\\
&=o\left(\frac{d_0}{\l _1}\right)+O\left( \frac{\e _{12}^2\left(Log(\e
  _{12}^{-1})\right)^{\frac{2(n-2)}{n}}}{\l _1d_0}\right)\nonumber\\
&=o\left(\frac{d_0}{\l _1} + \frac{\e _{12}}{\l _1d_0}\right).
\end{align}
Using \eqref{e:49},\eqref{e:410},\eqref{e:411} and \eqref{e:412},we
obtain
\begin{align*}
-\left(\n J(u), W_3\right) & \geq C\left(\frac{d_0}{\l
 _1}+\frac{d_0}{\l _2}+\frac{\e _{12}}{\l _1d_0}\right)\\
&\geq C\left(\frac{d_0}{\l _1}+\frac{d_0}{\l _2}+ \e
 _{12}^{(n-\frac{1}{2})/(n-2)}\right)\\
&\left(\sum\frac{1}{\l
    _i^2}+ \sum (1-J(u)^{\frac{n}{n-2}}a
  _i^{\frac{4}{n-2}}K(a_i)\right) + C\e _{12}^{\frac{n-\frac{1}{2}}{n-2}}.
\end{align*}
{\bf $3^{rd}$ region.} \, $\exists i\in \{1,2\}$ such that
$|1-J(u)^{\frac{n}{n-2}}a _i^{\frac{4}{n-2}}K(a_i)|\geq \frac{C'}{\l
  _i}.$\\
In this case, we set
$$
Z=-sign
\left(1-J(u)^{\frac{n}{n-2}}a_i^{\frac{4}{n-2}}K(a_i)\right)\tilde{\d}
_i.
$$
Using Proposition \ref{p:35}, we obtain
\begin{align*}
-\left(\n J(u), Z\right)&\geq
C|1-J(u)^{\frac{n}{n-2}}a_i^{\frac{4}{n-2}}K(a_i)| +
O\left(\frac{1}{\l _i}\right)+O(\e _{12})\\
&\geq \frac{C}{2}|1-J(u)^{\frac{n}{n-2}}a_i^{\frac{4}{n-2}}K(a_i)| +
\frac{C'}{4\l _i}+O(\e _{12}).
\end{align*}
We also set
$$
Z_1=
\begin{cases}
W_1 \qquad \mbox{if} \qquad M\l _1 \leq \l _2 \\
-\l _2 \frac{\partial\tilde\d _2}{\partial\l _2} \qquad
\mbox{if}\qquad 2M\l _1 \geq \l _2.
\end{cases}
$$
Setting $W_4 = Z + Z_1 \sqrt{C'}$, we derive
\begin{align*}
-\left(\n J(u), W_4\right)&\geq
C\left(|1-J(u)^{\frac{n}{n-2}}a_i^{\frac{4}{n-2}}K(a_i)| +
\sum \frac{1}{\l _k}+\e _{12}\right)\\
&\geq \gamma \left(\sum |1-J(u)^{\frac{n}{n-2}}a_k^{\frac{4}{n-2}}K(a_k)|^2 +
\sum \frac{1}{\l _k^2}+\e _{12}^{\frac{n-\frac{1}{2}}{n-2}}\right).
\end{align*}
Hence our global vector field will be built using a convex combination
of $W_1$, $W_2$, $W_3$ and $W_4$ and will satisfy obviously i., iii. and
iv. Next, we give the proof of ii. As in \cite{B1} and \cite{BCCH}, it
is easy to prove that
\begin{align*}
-\left(\n J(u+\bar{v}), \frac{\partial{\bar{v}}}{\partial (\a _i, a_i,
 \l _i)}(W)\right)&\leq c|\bar{v}| |\n J(u+\bar{v})|\\
&=O\left(|\bar{v}|^2 + |\n J(u+\bar{v})|^2\right).
\end{align*}
By Propositions \ref{p:33} and \ref{p:34}, we can prove that
$$
|\n J(u+\bar{v})|= O\left(\sum |1-J(u)^{\frac{n}{n-2}}a_k^{\frac{4}{n-2}}K(a_k)| +\sum \frac{1}{\l _k}+\e _{12} +|\bar{v}|\right).
$$
Using now the estimate of $\bar{v}$, we easily derive ii. Thus our proposition follows.
\hfill$\Box$\\
Next, we state the following result whose proof is similar to
corresponding statements in  Lemma A1.1 \cite{B1}
\begin{lemma}\label{l:B}(Lemma A1.1 \cite{B1})
Let $u=\a _{i_0}\d _{(x_{i_0}, \l _{i_0})}+\a _{j_0}\d _{(x_{j_0}, \l
  _{j_0})}$, where \\
(i)\quad $\e _{i_0 j_0} \geq \d _1$, $\d _1$ a given constant\\
(ii)\quad $B^2 \geq \l _{i_0}$, $\l _{j_0} \geq B$\\
If, given $\d _1$, $B$ is large enough, there is a pseudogradient
  vector field of $J$, built with the Yamabe gradient on $u$, which
  leads functions such as $u$ in the neighborhood of functions of the
  type $\a \d _{(y,\l )}+v$, where $y$ is close to
  $\frac{1}{2}(x_{i_0}+x_{j_0})$ up to $O(\frac{1}{B})$, $\l \geq cB$
  ($c$ is a universal constant) and $||v||=o(1)$.
\end{lemma}
Now, we will use the above
Lemma in the proof of the next main result.
\begin{proposition}\label{p:62}
Let $\e _0 > 0$ small enough. There exists a vector field $Z_0$
defined in
$$
W_{\e _0}=\{\a _1 \tilde\d _{x,\l _1} +\a _2 \tilde\d _{y_0, \l _2}/
\a _i \geq 0, \a _1+\a _2 =1, x\in X, \l _i > \e _0^{-1}, \e _{12}< \e
_0\}
$$
which can be extended to
$$
W(2,\e _0)=\{\a _1 \tilde\d _{a_1,\l _1} +\a _2 \tilde\d _{a_2, \l _2}
+ v \in V_b(2,\e _0)/ a_1, a_2 \in \overline{S^n_+}\}
$$
so that the following holds:\\
$f_\l (C_{y_0}(X))$ retracts by deformation on $X\cup W_u(y_0,
y_{i_0})_\infty \cup D$, where $D \subset \sigma $ is a stratified set
(in the topological sense, that is, $D\in \Sigma_k(S^n_+)$, the group of
chains of dimensions k) and where $\sigma = \cup _{y_i \in X\diagdown
  \{y_{i_0}, y_0\}}W_u(y_0, y_i)_\infty$ is a manifold in dimension at most
$k-1$.\\
Here $W_u$ denotes the unstable manifold for $Z_0$.
\end{proposition}
\proof
First, we notice that assumption $(H_2)$ implies that any critical
point $y$ of $K_1$ such that $y\in X$ satisfies $(\partial
K/\partial\nu ) > 0$. Now, we distinguish five cases\\
{\bf Case 1.} \quad There exists $i$ such that $\a _i$ is far away
from $J(u)^{\frac{-n}{4}}K(a_i)^{\frac{2-n}{4}}$($i=1,2$).\\
We set
$$
Z_1=\tilde\d _{a_i,\l _i} \dot{\a _i}\quad \mbox{with } \dot{\a _i}
=\begin{cases}
1 \quad \mbox{if } \a _i >
J(u)^{\frac{-n}{4}}K(a_i)^{\frac{2-n}{4}}+\eta\\
-1  \quad \mbox{if } \a _i <
J(u)^{\frac{-n}{4}}K(a_i)^{\frac{2-n}{4}}-\eta
\end{cases}
$$
where $\eta$ is a positive constant.\\
Using Proposition \ref{p:35}, we obtain
\begin{align*}
-\left(\n J(u), Z_1\right)&= c + O\left(\l _i^{-1}\right)+ O(\e _{12})\\ 
&\geq C > 0.
\end{align*}
{\bf Case 2.} \quad For each $i \in \{1,2\}$,  $\a
_i=J(u)^{\frac{-n}{4}}K(a_i)^{\frac{2-n}{4}}$ and $x\notin
\mathcal{V}(y_i,\rho)$, where \\
$\rho < \frac{1}{2}\min _{i\not=
  j}d(y_i,y_j)$
and $y_i$ is any critical point of $K_1= K_{/\partial S^n_+}$.\\
In this case, we have $d(x,y_0) \geq c$, thus $\e
_{12}=o\left(\frac{1}{\l _i}\right)$ for $i=1, 2$.\\
Two subcases may occur\\
If $\l _1 \leq C_1 \l _2$, where $C_1$ is a  large
enough positive constant, we set 
$$
Z_{21}= \frac{1}{\l _1}\frac{\partial\tilde\d _1}{\partial a_1} \n _TK.
$$
If $\l _1 \geq C_1 \l _2$, we set
$$
Z_{22}= Z_{21}+ \l _2 \frac{\partial\tilde\d _2}{\partial\l _2}.
$$
Using Propositions \ref{p:33} and \ref{p:34}, we derive
$$
-\left(\n J(u), Z_{2i}\right) \geq c \l _1^{-1} + c \l _2^{-1} +
\e _{12}\qquad \mbox{for }
i= 1,  2.
$$
{\bf Case 3.} \quad For each $i \in \{1,2\}$,  $\a
_i=J(u)^{\frac{-n}{4}}K(a_i)^{\frac{2-n}{4}}$ and $x\in
\mathcal{V}(y_i,2\rho)$, where \\
 $y_i$ is any critical point of $K_1$ such that $y_i\not= y_0$.\\
Since $x\in X$, $y_i \in X$ and therefore $(\partial K/\partial\nu
)(y_i) > 0$. Now, we set
$$
Z_3 = \l _1 \frac{\partial\tilde\d _1}{\partial\l _1} +\l _2
\frac{\partial\tilde\d _2}{\partial\l _2}.
$$
Using Proposition \ref{p:33}, we obtain
$$
-\left(\n J(u), Z_{3}\right) \geq c \l _1^{-1}+ c \l _2^{-1}+ \e _{12}. 
$$
{\bf Case 4.} \quad For each $i \in \{1,2\}$,  $\a
_i=J(u)^{\frac{-n}{4}}K(a_i)^{\frac{2-n}{4}}$ and $x\in
\mathcal{V}(y_0,2\rho)$\\
In this case, we use the vector field defined in the proof of
Proposition \ref{p:61} which we combine with the vector field
defined in Lemma \ref{l:B}.\\
{\bf Case 5.} \quad $\a _1 = 0$ or $\a _2 = 0.$\\
In this case, we only have one mass and we use the vector field
defined in Lemma \ref{l:43}.\\ \\
Our global vector field $Z_0$ will be built using a convex combination
of vector fields defined in cases 1- 5.\\
Now, let $u= \a \d _{x,\l} + (1-\a )\d _{y_0,\l} \in f_\l
(C_{y_0}(X))$.\\
The action of the flow of the pseudogradient $Z_0$ is described as
follow.\\
If $\a < 1/2$, the flow of $Z_0$ brings $\a$ to zero, and thus in this
case $u$ goes to $\overline{W_u((y_0)_\infty)} = \{y_0\}$.\\
If $\a > 1/2$, the flow of $Z_0$ brings $\a$ to $1$, and thus $u$ goes,
in this case, to $\overline{W_u((y_{i_0})_\infty)} =X$.\\
If $\a = (1-\a ) = 1/2$, we have an action on $x\in X=
\overline{W_s(y_{i_0})}$. In this case, $u$ goes to $W_s(y_i)$, where $y_i$
is a critical point of $K_1$ dominated by $y_{i_0}$ and two cases may
occur :\\
In the first case $y_i \not= y_0$, then $x$ goes to $W_u(y_0,
y_i)_\infty$.\\
In the second case $y_i=y_0$, $u$ goes to $W_u(y_0)_\infty$ by the
vector field defined in Lemma \ref{l:B}.\\
Then our result follows. 
\hfill$\Box$\\ \\
We now prove our theorems.

\begin{pfn}{\bf Theorem \ref{t:13}}
We argue by contradiction. Assume that (1) has no solution. The strong
retract defined in Proposition \ref{p:62} does not intersect
$W_u(y_0,y_{i_0}))_\infty$ and thus it is contained in $X\cup D$ (see
Proposition \ref{p:62}). Therefore $H_*(X\cup D) = 0$, for all $*\in
\N^*$, since $f_\l (C_{y_0}(X))$ is a contractible set.\\
Using the exact homology sequence of $(X\cup D, X)$, we have 
$$
...\to H_{k+1}(X\cup D) \to ^{\pi}  H_{k+1}(X\cup D, D) \to
^{\partial}  H_k(X) \to ^{i}  H_k(X\cup D) \to ...
$$
Since $H_*(X\cup D) = 0$, for all $*\in \N^*$, then
$H_k(X)=H_{k+1}(X\cup D, X)$.\\
In addition, $(X\cup D, X)$ is a stratified set of dimension at most
$k$, then $H_{k+1}(X\cup D,X) = 0$, and therefore $H_k(X)=0$. This
yields a contradiction since $X$ is a manifold in dimension $k$
without boundary. Then our theorem follows.
\end{pfn}\\
\begin{pfn}{\bf Theorem \ref{t:14}}
Assume that (1) has no solution. By the above arguments, $X\cup \left(\cup
_{y_i\in B_k}W_u(y_0, y_i)\right)\cup D$ is a strong retract of $f_\l
(C_{y_0}(X))$, where $D\subset \sigma$ is a stratified set and where
$\sigma = \cup _{y_i\in X\diagdown \left(B_k \cup \{y_{i_0}\}\right)}W_u(y_0, y_i)_\infty$
is a manifold in dimension at most $k$.\\
Since $\mu _i (y_0)=0$ for each $y_i \in B_k$, $f_\l (C_{y_0}(X))$
retracts by deformation on $X\cup D$, and therefore $H_*(X\cup D)=0$,
for all $*\in \N^*$. Using the exact homology sequence of $(X\cup D,
D)$, we obtain $H_{k+1}(X\cup D, X) = H_k(X) = 0$, a contradiction,
and therefore our result follows.
\end{pfn}


\begin{thebibliography}{99}

\bibitem{ALM} A. Ambrosetti , Y.Y. Li , A. Malchiodi, 
\emph{Yamabe and Scalar
Curvature Problems under boundary conditions}, preprint S.I.S.S.A.,
ref. 52/2000/M. Preliminary note on C.R.A.S., S{\'e}rie 1 330 (2000),
1013-1018.
\bibitem{AB}
 T. Aubin and A. Bahri,
\emph{ Methodes de topologie algebrique pour le probl{\`e}me de 
la courbure scalaire prescrite,} J. Math. Pures et Appl.  \textbf{76} (1997), 525--549.
\bibitem{B1}
 A. Bahri,
\emph{An invarient for Yamabe-type flows with applications to scalar curvature
problems in high dimension,} A celebration of J. F. Nash Jr., Duke Math. J. 
\textbf{81} (1996), 323-466.
\bibitem{B2}
 A. Bahri,
\emph{ Critical point at infinity in some variational problems,}
 Pitman Res. Notes Math, Ser \textbf{ 182}, Longman Sci. Tech. Harlow 1989.
\bibitem{BC1} 
 A. Bahri and J. M. Coron,
\emph{The scalar curvature problem on the standard three dimensional spheres,}
 J. Funct. Anal. \textbf{95} (1991), 106-172.
\bibitem{BC2}
A. Bahri  and J.M. Coron,
\emph{ On a nonlinear elliptic equation involving the critical Sobolev exponent : the effect of topology of the domain ,}
 Comm. Pure Appl. Math.\textbf{ 41}(1988), 255--294.
\bibitem{BR}
A. Bahri and P. Rabinowitz,
\emph{Periodic orbits of hamiltonian systems of three body type, }
Ann. Inst. H. Poincar{\'e} Anal. Non lin{\'e}aire \textbf{8} (1991), 561-649.
\bibitem{BCCH}
 M. Ben Ayed, Y. Chen, H. Chtioui and M. Hammami,
\emph{ On the prescribed scalar curvature problem on 4-manifolds,}
 Duke Math. J. \textbf{84} (1996), 633-677.
\bibitem{BCH}
 Ben Ayed, H. Chtioui and M. Hammami,
\emph{The scalar curvature problem on higher dimensional spheres,}
 Duke Math. J. \textbf{93} (1998), 379-424.
\bibitem{BP}
G. Bianchi and X. B. Pan,
\emph{Yamabe equations on half spheres,}
Nonlinear Anal. \textbf{37} (1999), 161-186.
\bibitem{BrC}
H. Brezis and J. M. Coron,
\emph{Convergence of solutions of H-systems or how to blow bubbles,}
 Arch. Rational Mech. Anal. \textbf{89} (1985), 21-56.
\bibitem{CGY}
S. A. Chang, M. J. Gursky and P. Yang,
\emph{The scalar curvature equation on 2 and 3 spheres,}
 Calc. Var. \textbf{1} (1993), 205-229.
\bibitem{CY}
S. A. Chang and P. Yang,
\emph{A perturbation result in prescribing scalar curvature on $S^n$,}
Duke Math. J. \textbf{64} (1991), 27-69.
\bibitem{C}
P. Cherrier,
\emph{Probl{\`e}mes de Neumann non lin{\'e}aires sur les vari{\'e}t{\'e}s 
Riemaniennes,}
J. Funct. Anal. \textbf{57} (1984), 154-207.
\bibitem{DMA}
Z. Djadli, A. Malchiodi and M. Ould Ahmedou,
\emph{Prescribing the scalar and the boundary mean curvature on the three 
dimensional half sphere,} 
Preprint \textbf{} (2001).
\bibitem{E1}
J. Escobar,
\emph{Conformal deformation of Riemannian metric to scalar flat metric
 with constant  mean curvature on the boundary,}
Ann. of Math. \textbf{136} (1992), 1-50.
\bibitem{E2}
J. Escobar,
\emph{Conformal metrics with prescribed mean curvature on the boundary,}
Cal. Var. \textbf{4} (1996), 559-592.
\bibitem{HL1}
Z. C. Han and Y.Y. Li,
\emph{The Yamabe problem on manifolds with boundaries : existence and 
compactness results,}
 Duke Math. J. \textbf{99} (1999), 489-542.
\bibitem{HL2}
Z. C. Han and Y.Y. Li,
\emph{The existence of conformal metrics with constant scalar curvature 
and constant boundary mean curvature,}
Comm. Anal. Geom. \textbf{8} (2000), 809-869.
\bibitem{H}
E. Hebey,
\emph{Changements de metriques conformes sur la sphere, le probl{\`e}me 
de Nirenberg,}
Bull. Sci. Math. \textbf{114} (1990), 215-242.
\bibitem{L1}
Y.Y. Li,
\emph{The Nirenberg problem in a domain with boundary,}
Top. Meth. Nonlin. Anal. \textbf{6} (1995), 309-329.
\bibitem{L2}
Y.Y. Li,
\emph{Prescribing scalar curvature on $S^n$ and related topics,
Part I,} J. Diff. Eq. \textbf{120} (1995), 319-410.
\bibitem{L3}
Y.Y. Li,
\emph{Prescribing scalar curvature on $S^n$ and related topics, Part II :
 existence and compactness,}
Comm. Pure Appl. Math. \textbf{49} (1996), 437-477.
\bibitem{L}
P. L. Lions,
\emph{The concentration compactness principle in the calculus of variations.
 The limit case,}
Rev. Mat. Iberoamericana \textbf{1} (1985), I: 165-201; II: 45-121.
\bibitem{M}
J. Milnor,
\emph{Lecturess on h-Cobordism Theorem,}
 Princeton University Press, Princeton \textbf{} 1965.
\bibitem{R}
O. Rey,
\emph{Boundary effect for an elliptic Neumann problem with critical 
nonlinearity,}
Comm. Partial Diff. Eq. \textbf{22} (1997), 1055-1139.
\bibitem{SZ}
R. Schoen and D. Zhang,
\emph{Prescribed scalar curvature on the n-sphere,}
Calculus of Variations and Partial Differential Equations, \textbf{4} (1996),
1-25.
\bibitem{S}
M. Struwe,
\emph{A global compactness result for elliptic boundary value problems 
involving nonlinearities,}
Math. Z. \textbf{187} (1984), 511-517. 
\end{thebibliography}
\end{document}